\documentclass{amsart}
\usepackage{amssymb,amsmath,amscd}

\usepackage[colorlinks=true]{hyperref}

\newtheorem{thm}{Theorem}[section]

\newtheorem{corollary}[thm]{Corollary}
\newtheorem{prop}[thm]{Proposition}
\newtheorem{lemma}[thm]{Lemma}
\newtheorem{fact}[thm]{Fact}

\theoremstyle{definition}
\newtheorem{defn}[thm]{Definition}
\newtheorem{example}[thm]{Example}

\theoremstyle{remark}
\newtheorem{remark}[thm]{Remark}

\newcommand{\bt}{\begin{thm}}
\newcommand{\et}{\end{thm}}
\newcommand{\bp}{\begin{prop}}
\newcommand{\ep}{\end{prop}}
\newcommand{\bd}{\begin{defn}}
\newcommand{\ed}{\end{defn}}
\newcommand{\bl}{\begin{lemma}}
\newcommand{\el}{\end{lemma}}
\newcommand{\bfa}{\begin{fact}}
\newcommand{\efa}{\end{fact}}
\newcommand{\bc}{\begin{corollary}}
\newcommand{\ec}{\end{corollary}}
\newcommand{\bex}{\begin{example}}
\newcommand{\eex}{\end{example}}
\newcommand{\br}{\begin{remark}}
\newcommand{\er}{\end{remark}}
\newcommand{\ben}{\begin{enumerate}}
\newcommand{\een}{\end{enumerate}}

\newcommand{\sotto}[2]{#1_{#2}}

\newcommand{\rrr}{\rightarrow}
\newcommand{\ra}{\rightarrow}

\newcommand{\ds}{\displaystyle}

\newcommand{\ideal}[1]{\sotto {{\mathcal I}}{#1}}

\newcommand{\exact}[3]
{0 \rrr #1 \rrr #2
\rrr #3 \rrr 0}

\newcommand{\PP}{\mathbb{P}}

\newcommand{\Z}{\mathbb{Z}}

\newcommand{\coo}{{\mathcal O}}

\newcommand{\K}{\mathbb{K}}

\begin{document}

\title[The Hilbert schemes of curves may after all be connected]{The Hilbert schemes of locally Cohen-Macaulay curves in $\PP^3$
may after all be connected}

\author[P.~Lella]{Paolo Lella}
\address{Dipartimento di Matematica dell'Universit\`{a} di Torino\\
         Via Carlo Alberto 10,
         10123 Torino, Italy}
\email{\href{mailto:paolo.lella@unito.it}{paolo.lella@unito.it}}
\urladdr{\href{http://www.personalweb.unito.it/paolo.lella/}{www.personalweb.unito.it/paolo.lella/}}

\author[E.~Schlesinger]{Enrico Schlesinger}
\address{Dipartimento di Matematica del Politecnico di Milano\\
         Piazza Leonardo da Vinci 32,
         20133 Milano, Italy}
\email{\href{mailto:enrico.schlesinger@polimi.it}{enrico.schlesinger@polimi.it}}

\thanks{
The second author was partially supported by
MIUR PRIN 2007:
  {\em 	 Moduli, strutture geometriche e loro applicazioni}.
 }


\begin{abstract}
Progress on the problem whether the Hilbert schemes of locally Cohen-Macaulay curves in $\PP^3$
are connected has been hampered by the lack of an answer to a question raised by Robin Hartshorne
in \cite{hconn} and more recently in  \cite{aim}: does there exist a flat irreducible family of curves
whose general member is a union of $d \geq 4$ disjoint lines on a smooth quadric surface and
whose special member is a locally Cohen-Macaulay curve in a double plane?
In this paper we give a positive answer to this question: for every $d$ we construct a family with the required properties,
whose special
fiber is an {\em extremal} curve in the sense of \cite{extremal}.
From this we conclude that every effective divisor in a smooth quadric surface 
is in the connected component of its Hilbert scheme that contains extremal curves.
\end{abstract}

\maketitle

\section{Introduction}
Liaison theory is a major tool in the study of Hilbert schemes of codimension two subschemes
of a fixed projective scheme. We refer to  \cite{htoul} and \cite{migliore} for a modern treatment
of this theory. This paper is concerned with the case of curves in $\PP^3$.  By the term {\em curve}
we will mean a one dimensional subscheme without isolated or embedded zero dimensional components;
in the literature such an object is usually called a {\em locally Cohen-Macaulay}  curve.
Let us recall the main results of biliaison theory in the context of curves in $\PP^3$.
We say that a curve $D$ is obtained from a curve $C$ by {\em elementary biliaison}
of {\em height} $h$ if $D$ is linearly equivalent to $C+hH$ on a surface
$S$ that contains both $C$ and $D$ and has $H$ as its plane section \cite{hgd} (see also \cite[Chapter III]{MDP}).
For example, an effective divisor $D$ of bidegree $(a,b)$ on a smooth quadric surface $Q \cong \PP^1 \times \PP^1 \subset \PP^3$ is obtained
by an elementary biliaison of nonnegative height from a curve $C$ which is the disjoint union of $d=|a-b|$ lines on
the quadric.
Biliaison is the equivalence relation generated by elementary biliaisons. Rao \cite{rao} proved that there
is an invariant that distinguishes
biliaison equivalence classes: the finite length graded module
$M_C = \bigoplus_{n \in \Z} H^1 (\PP^3, \ideal{C} (n))$, which is now commonly referred to as
the Rao (or Hartshorne-Rao) module of  $C$; two curves are in the same biliaison
class if and only if their Rao modules are isomorphic up to a twist. The structure of a biliaison
class is well understood. A curve $C$ in a biliaison class is said to be {\em minimal} if for every other
curve $D$ of the class $M_D \cong M_C (-h)$ with $h \geq 0$; this means the Rao module of $D$ is obtained shifting
the Rao module of $C$ to the right. The main result of the theory is:
\begin{enumerate}
\item[i)]
Every biliaison class contains minimal curves; the family of the minimal curves of the class is irreducible, and
any two minimal curves of the class can be joined by a finite number of elementary biliaisons of height zero.
\item[ii)]
If $D$ is a non minimal curve, then $D$ is obtained from a minimal curve $C$ of the class by a finite
sequence of elementary biliaisons of positive height. This is known as the {\em Lazarsfeld-Rao property}.
\end{enumerate}
Note that from $ii)$ it follows that a curve is minimal in its biliaison class if and only if it has minimal degree among curves
of the class. This result has a long history. Lazarsfeld and Rao \cite{lazrao} proved the Lazarsfeld Rao property under a cohomological
condition on $C$ that guarantees $C$ is minimal; they only considered {\em trivial} elementary biliaisons ($C_i=C_{i-1}+nH$, no linear
equivalence allowed), but at the end of the process they needed a deformation with constant cohomology.
The existence of minimal curves and the Lazarsfeld-Rao property were proven independently in \cite{MDP}
and \cite{BBM}. Strano \cite{strano} showed the deformation at the end of the Lazarsfeld-Rao process is not needed if
one allows linear equivalence in the definition of elementary biliaison.  The version of the theorem we have given is due to Hartshorne
\cite[Theorem 3.4]{htoul}, where the precise conditions on the ambient projective scheme are determined for the Lazarsfeld-Rao
property to hold for biliaison classes of codimension $2$ subschemes.

We would like to stress the fact that for this theory to work it is necessary to consider (locally Cohen-Macaulay) curves:
even if one starts with a smooth irreducible curve, the minimal curve of the class may fail to be reduced or irreducible;
it may even not be generically a  local complete intersection (as it was assumed by Rao in \cite{rao}). So
let $H_{d,g}$ denote the Hilbert scheme parametrizing (locally Cohen-Macaulay) curves of degree $d$ and arithmetic genus $g$ in $\PP^3$; it is an open
subscheme of the full Hilbert scheme. In  \cite{bounds} and \cite{extremal} Martin-Deschamps and Perrin have shown that, whenever nonempty,
$H_{d,g}$ has an irreducible component $E_{d,g}$ whose closed points correspond to curves that have maximal cohomology, in the sense that
\begin{equation}\label{extr-in}
\dim H^{i} (\PP^3,\ideal{C}(n)) \leqslant \dim H^{i} (\PP^3,\ideal{E}(n))
\end{equation}
for every $C \in H_{d,g}$, $E \in E_{d,g}$, $n \in \Z$ and $i=0,1,2$. Note that by (\ref{extr-in}) there is no obstruction
coming from the semicontinuity of cohomology to specializing any curve in $H_{d,g}$ to an extremal curve. This remark
raised the question whether every curve can be specialized to an extremal curve. This is known to be false
for curves that are not generically of embedding dimension two \cite{ns-comp}, but it is open for, say, smooth curves.
A weaker version of this question,
proposed by Hartshorne in \cite{hconn} and \cite{hconn2},  is whether $H_{d,g}$ is connected, that is, if every curve in $H_{d,g}$ belongs to the connected component containing the extremal curves.
It is an interesting problem because the Hilbert scheme of smooth curves is not connected \cite{AG}, while the full Hilbert scheme
is connected \cite{hthesis}, but through schemes with zero dimensional components.

Let us review what is known about the problem of connectedness of $H_{d,g}$. No example of a nonconnected $H_{d,g}$
has been found so far. The Hilbert
scheme $H_{d,g}$ is connected when $g \geqslant \binom{d-3}{2}-1$  (see \cite{AA1,hmdp3,sabadini})
and when $d \leqslant 4$ (see \cite{nthree,ns-comp};
this is non trivial, because when $g$ is very negative the Hilbert scheme has a large number of irreducible components).
Hartshorne \cite{hconn} has shown that smooth irreducible nonspecial curves
and arithmetically Cohen-Macaulay (ACM) curves are in the connected component containing the extremal curves,
and so are Koszul curves by \cite{perrin-pas}. Building on the results of \cite{hconn}, the second author has shown
\cite{footnote} that, if a minimal curve $C$ can be connected to an extremal curve by flat families
lying on surfaces of degree $s$, where $s$ is the least degree of a surface containing $C$, then {\em every} curve in the biliaison class of $C$
is in the connected component of  extremal curves in its Hilbert scheme.

By \cite{hgd} and \cite{2h}, any curve contained in a singular quadric surface, including a double plane, is in the connected component of  extremal curves \cite{2h}. On the other hand, the case of curves on a smooth quadric surface
has been open so far. By biliaison it is enough to deal with divisors of bidegree $(d,0)$. The case $d \leq 2$ is then trivial,
and Nollet \cite{nthree} has shown that is possible  to specialize a divisor of bidegree $(3,0)$ to an extremal curve, but
it has been an open question whether one could specialize four (or more) disjoint lines on a smooth quadric
to an extremal curve of the same genus; consequently, the Hilbert scheme of curves of degree $H_{10,12}$, which has an irreducible component
whose general member is a divisor of bidegree $(7,3)$ on a smooth quadric surface, was proposed as a  candidate for an example
of a nonconnected $H_{d,g}$: see \cite[Ex. 4.2]{hconn}, \cite[Section 4]{hconn2}, and the open problems list
of the 2010 AIM workshop on \emph{Components of  Hilbert Schemes}  available at \href{http://aimpl.org/hilbertschemes}{http://aimpl.org/hilbertschemes}.

In this paper we give a positive answer to the above question:
\bt
For every $d \geq 3$, there exists a flat family of curves over the affine line
whose general member $C$ is a union of $d$ disjoint lines on a smooth quadric surface and
whose special member $E$ is an extremal curve in the double plane. It follows that every curve in the biliaison
class of $C$, in particular every curve on a smooth quadric surface,
is in the connected component of the extremal curves in its Hilbert scheme.
\et
The paper is structured as follows. In section \ref{two}, where we fix notation and terminology,
we observe that there are degree $d$ extremal curves, supported on a line, whose defining equations
are homogeneous with respect to the weight vector $\omega=(d,2,1,1)$. In section \ref{three} we consider
curves $C$ that are the union of $d$ disjoint lines on a smooth quadric surface, and
we determine conditions on $C$ under which the initial ideal $\text{in}_\omega(I_C)$ is
(after saturation) the ideal of an extremal multiple line $E$.
It is then a standard fact that there is a flat family having $C$ as generic fiber  and $E$ as special fiber.
One should note that, while the computation of the initial ideal is specific to the case under consideration,
the weight $\omega$ only depends on the structure of extremal curves, and we will show in a future paper
that our procedure can be applied in a much more general setting. 

\vspace{0.5cm} \noindent {\bf Acknowledgments.}
We thank Robin Hartshorne for many suggestions that allowed us to improve the first draft of the paper.
The second author would like to thank the American Institute of Mathematics  and the organizers of the AIM workshop {\em Components of Hilbert Schemes} Robin Hartshorne, Diane Maclagan, and Gregory Smith.  The seeds of this paper were planted through fruitful discussions with the other participants.

\section{Extremal curves} \label{two}
In this section we establish notation and terminology and review some known results that we will need later.
We work over an algebraically closed field $\K$ of arbitrary characteristic. We denote with the symbol
$\ideal{X}$ the ideal sheaf of a subscheme $X \subset \PP^3$. Given a coherent
sheaf $\mathcal{F}$ on $\PP^3$, we define $h^i(\mathcal{F})= \dim \, H^{i}(\PP^3,\mathcal{F})$ and
$\ds H^i_* (\mathcal{F})= \bigoplus_{n \in \Z} H^{i}(\PP^3,\mathcal{F}(n))$. We write
$R=\K[x,y,z,w]$ for the homogeneous coordinate ring $H^{0}_* (\coo_{\PP^3})$ of $\PP^3$.

\bd
A {\em curve} in $\PP^3$ (or more precisely a locally Cohen-Macaulay curve)
is a one dimensional subscheme $C \subset \PP^3$ without zero dimensional associated points; this means
that all irreducible components of $C$ have dimension $1$, and that $C$ has no embedded points.
\ed
We denote by $H_{d,g}$ the Hilbert scheme parametrizing curves of degree $d$ and arithmetic
genus $g$ in $\PP^3$ \cite{MDP}. This is an open subscheme of the full Hilbert scheme parametrizing
all one dimensional subschemes of $\PP^3$ with Hilbert polynomial $dn+1-g$.

\bd We say that a curve $E \subset \PP^3$ of degree $d$ and genus $g$ is {\em extremal} if
\begin{equation}\label{extr_def}
h^{i} (\ideal{C}(n)) \leqslant h^{i} (\ideal{E}(n))
\end{equation}
for $i=0,1,2$, for every curve $C$ of degree $d$ and genus $g$, and for every $ n \in \Z$,
\ed
Thus a curve is extremal if it has the largest possible cohomology.
One knows that $H_{d,g}$ is nonempty if  and only if either $g=\binom{d-1}{2}$, in which case it is
irreducible and consists of plane curves, or $g \leqslant \binom{d-2}{2}$. Whenever nonempty,
$H_{d,g}$ contains extremal curves \cite{bounds}; in fact, the extremal curves form an irreducible
$E_{d,g}$ of $H_{d,g}$ \cite{extremal}.

\br
Our definition of extremal curves is equivalent to the one given by Hartshorne \cite{hconn,hconn2}.
Martin-Deschamps and Perrin's definition of an extremal curve is more restrictive, because they
require an extremal curve not to be arithmetically Cohen-Macaulay (ACM);
the difference comes up only when $g=\binom{d-1}{2}$ and $g=\binom{d-2}{2}$, the only two cases
in which $H_{d,g}$ consists entirely of ACM curves.
As pointed out by Nollet \cite[Proposition 2.1]{subextremal},  the function $h^{1} (\ideal{E} (n))$ for an extremal curve $E$
can be read from  \cite{bounds}, where the
bounds (\ref{extr_def}) are proven  under  the assumption that the field $\K$ has characteristic zero;
this assumption on the characteristic is not necessary \cite{subextremal}.
\er

Martin-Deschamps and Perrin also
compute the Rao module $M_E = H^{1}_* (\ideal{E})$ of an extremal curve:

\bt[\cite{bounds,extremal}]
Let $(d,g)$ be two integers satisfying $d \geqslant 2$ and
$g \leqslant \binom{d-2}{2}-1$.
A curve $E$ of degree $d$ and genus $g$ is extremal if and only if
$$
M_E \cong R/(l_1,l_2,F,G) (b)
$$
where $(l_1,l_2,F,G)$ is a regular sequence, $\deg(l_1)=\deg(l_2)=1$,
$\deg(F)= \binom{d-2}{2}-g$, $\deg(G)= \binom{d-1}{2}-g$, and
$b=\deg(F)-1$.
\et

The following proposition describes extremal curves supported on a line. It is a special case of
\cite[Proposition 0.6]{extremal};  we state it in the form
needed later in the paper, and we prove it for convenience of the reader.

\bp \label{extr_mls}
Let $(d,g)$ be a pair of integers satisfying $g \leqslant \binom{d-2}{2}-1$. Let $F$ and $G$
be two forms of degrees $\deg(F)= \binom{d-2}{2}-g$ and $\deg(G)= \binom{d-1}{2}-g$ in $\K[z,w]$ with no common zeros.
The surface $S$ of equation $xG-y^{d-1}F=0$ is irreducible and generically smooth along the line
$L$ of equations $x=y=0$. It therefore contains a unique curve $E$ of degree $d$ supported on $L$.
The curve $E$ is extremal of degree $d$ and genus $g$, and its Rao module is
$$
M_E \cong R/(x,y,F,G) (b) \cong \K[z,w]/(F,G) (b)
$$
where $b=\deg(F)-1$. The homogeneous ideal of $E$ is generated by $x^2$, $xy$, $y^d$ and $xG-y^{d-1}F$.
\ep

\begin{proof}
The surface $S$ is irreducible because $F$ and $G$ have no common zeros, and it is smooth at points
of $L$ where $G$ is different from zero. Therefore  the ideal of $L$ in the local ring
$\coo_{S,\xi}$ of the generic
point $\xi$ of $L$ is generated by one function $t$ , and the ideal of a curve of degree $d$ supported on $L$
must be $t^d\coo_{S,\xi}$. Since a locally Cohen-Macaulay curve supported on $L$ is determined by its ideal
at the generic point of $L$, we see that there is a unique curve $D_m \subset S$ supported on $L$ of degree $m$
for every $m \geqslant 1$. For $m=d-1$, the curve $P=D_{d-1}$ is the planar multiple structure of equations
$x=y^{d-1}=0$. 
We note that
$\ideal{P} \otimes \coo_L \cong \coo_L (-1) \oplus \coo_L (1-d)$  where the two generators are the images
of $x$ and $y^{d-1}$. The two forms $F$ and $G$ define a surjective map $\coo_L (-1) \oplus \coo_L (1-d) \ra \coo_L (b)$;
composing this with the natural map  $\ideal{P} \ra \ideal{P} \otimes \coo_L$  we obtain
a surjection $\phi:\ideal{P} \ra \coo_L (b)$. We let $E$ be the subscheme of $\PP^3$ whose ideal sheaf is
the kernel of $\phi$. By construction we have an exact sequence
 \begin{equation}\label{extr_seq}
    \exact{\ideal{E}}{\ideal{P}}{\coo_L (b)}
 \end{equation}
This sequence shows that $E$ is a (locally Cohen-Macaulay) curve of degree $d$ and genus $g$, and that its homogeneous
ideal is generated by $x^2$,$xy$,$y^d$ and $xG-y^{d-1}F$. Therefore $E=D_d$ is the unique curve of degree $d$ contained
in $S$ and supported on $L$. Finally, the long exact cohomology sequence associated to (\ref{extr_seq}) shows that the Rao
module of $E$ is
$$
M_E= \K[z,w] /(F,G) (b) = R/(x,y,F,G) (b).
$$
Hence $E$ is an extremal curve.
\end{proof}

\section{The main result} \label{three}
In this section we construct a flat family of curves whose general member is a disjoint union
of lines on a smooth quadric surface $Q$ and whose special member is an extremal multiple line.
This specialization is obtained considering the initial ideal with respect to the weight vector
$\omega=(d,2,1,1)$. In case $d=3$, such a specialization was constructed by Nollet \cite{nthree}
for a triple structure $3L$ on the line $L$ on $Q$ (without using the language of weight vectors and initial ideals).
However, for $d \geqslant 4$, it can be shown that this method fails for the $d$-uple structure $dL$ on $Q$,
because the initial ideal of $dL$ with respect to the weight vector $\omega=(d,2,1,1)$
defines a curve with embedded points. We now show that the method succeeds if one replaces $dL$ with
a sufficiently general divisor $C$ of bidegree $(d,0)$ on $Q$; the heart of the argument
is the proof that the initial ideal of $C$  with respect to the weight vector
$\omega=(d,2,1,1)$ defines an an extremal curve.


\bt \label{main}
Let $Q$ be the quadric surface of equation
$$
q(x,y,z,w)=x(x+w)-yz=0.
$$
For every $a \in \K$ let $L_a \subset Q$ denote
the line of equations
$x-az=y-a(az+w)=0$.
Given $d \geqslant 3$ and $a_1, \ldots, a_d \in \K$, consider the divisor
$$C=L_{a_1}+ \cdots + L_{a_d}$$
on $Q$.
If the sums $a_{i} + a_{j}$  for
$1 \leqslant i<j \leqslant d$ are all distinct, then there is a flat family of pairs $\mathcal{C} \subset \mathcal{Q} \ra \mathbb{A}^1$,
whose fiber over $1$ is $(C,Q)$, whose fiber over $t \neq 0$ consists of $d$ disjoint lines on a smooth
quadric surface, and whose fiber over $0$ is an extremal curve in the double plane of equation $x^2=0$.
\et

\begin{proof}
The basic remark is that a homogeneous polynomial $xG-y^{d-1}F$, where $F$ and $G$ are forms in $z$ and $w$,
is also homogeneous with respect to the weight vector $\omega=(d,2,1,1)$: indeed, if $\deg(F)=a$ and $\deg(G)=a+d-2$,
then
$$
\deg_\omega (xG)= d+a+d-2= 2d-2+a  =\deg_\omega (y^{d-1}F).
$$
Given a polynomial $P(x,y,z,w)$, the {\em initial form} $\text{in}_{\omega} (P)$ of $P$ with respect to $\omega$ is the sum of all the terms
$c_{\mathbf{n}} x^{n_1}y^{n_2}z^{n_3}w^{n_4}$ in $P$ for which the scalar product
$$
\omega \cdot \mathbf{n}= dn_1+2n_2+n_3+n_4
$$
is maximal. The initial ideal $\text{in}_\omega (I)$ of an ideal $I$ is the ideal generated by the initial forms  $\text{in}_{\omega} (P)$
as $P$ varies in $I$. There is a flat family over the affine line $\mathbb{A}^1$ whose fibers over $t \neq 0$
are isomorphic to $\textnormal{Proj}\, (R/I)$, while the special fiber over zero is the subscheme of $\PP^3$ defined
by the (saturation of)  $\text{in}_\omega (I)$: see for example \cite{bayer-m}
and \cite[Theorem 15.17]{eisenbud}. Roughly, this family is defined letting the one dimensional torus act on $\PP^3$
by $t[x:y:z:w]=[t^{d}x:t^2y:tz:tw]$ and taking the limit as  $t$ goes to zero,
so that set theoretically we are projecting $C$ onto the line $L$ of equations $x=y=0$,
but what is interesting is the scheme theoretic structure of the limit.

Let $C_0$ denote the subscheme defined by (the saturation of) $\text{in}_\omega(I_C)$,  where $I_C$ denotes the homogeneous ideal of
$C$. By the above remark there is a flat specialization from $C$ to $C_0$; since $\text{in}_{\omega} (q)=x^2$, the smooth
quadric $Q$ specializes to the double plane $x^2=0$ as $C$ specializes to $C_0$.
We let $l_i=x-a_iz$ and $m_i=y-a_i(a_iz+w)$ denote the given equations for the line $L_{a_i}$. Since
$I_C$ contains the product of the ideals of the lines $L_{a_i}$,
$$
y^d= \text{in}_{\omega} (m_1 \cdots m_d) \in I_{C_0}
$$
Therefore $C_0$ is contained in the complete intersection $x^2=y^d=0$ and so it is supported on the line $L$.
We will show that $I_C$ contains a polynomial $A(x,y,z,w)$ such that $\text{in}_\omega(A)=xG-y^{d-1}F$
where $F$ and $G$ are homogeneous forms in $\K[z,w]$, $\deg G= \binom{d}{2}$, $\deg F= \binom{d-1}{2}+1$,
and $F$ and $G$ have no common zero in $\PP^1=\textnormal{Proj}\, (\K[z,w])$. It follows that $C_0$ is contained
in the surface $S$ of equation $xG-y^{d-1}F$. By flatness, the Hilbert polynomial of $C_0$ coincides
with that of $C$, so $C_0$ is a one dimensional subscheme of $\PP^3$ of degree $d$ and genus $1-d$.
Let $E$ the largest Cohen-Macaulay curve contained in $C_0$: it is the curve of degree $d$ obtained
from $C_0$ throwing away its embedded points. By Proposition \ref{extr_mls} $E$ is the unique
curve of degree $d$ contained in $S$ and supported on the line $L$; it is an extremal curve of degree
$d$ and genus $1-d$. Since $E \subset C_0$ and the two schemes have the same Hilbert polynomial, we conclude
$E=C_0$. Thus the limit is an extremal curve, and the statement is proven.

To conclude the proof, we need to find $A \in I_C$ with $\text{in}_\omega(A)=xG-y^{d-1}F$. In principle,
this is a Groebner basis calculation, but luckily we can bypass such a calculation because we were
able to find, with the help of some computations performed with Macaulay 2 \cite{M2}, a
determinantal formula for $A$ (note that, while $I_C$ is generated by forms of degree $\leqslant d$,
the degree of $A$ is much larger than $d$). Let $A=xG-zB$ denote the determinant
\begin{equation}\label{ad}
A=\begin{vmatrix}
l_1         & m_1     & m_1^2  & \ldots & m_1^{d-1} \\
l_2         & m_2     & m_2^2  & \ldots & m_2^{d-1} \\
\vdots    & \vdots  & \vdots & \ddots & \vdots    \\
l_d         & m_d     & m_d^2  & \ldots & m_d^{d-1} \\
\end{vmatrix}
= x
\begin{vmatrix}
1         & m_1     & m_1^2  & \ldots & m_1^{d-1} \\
1         & m_2     & m_2^2  & \ldots & m_2^{d-1} \\
\vdots    & \vdots  & \vdots & \ddots & \vdots    \\
1         & m_d     & m_d^2  & \ldots & m_d^{d-1} \\
\end{vmatrix}
-z
\begin{vmatrix}
a_1         & m_1     & m_1^2  & \ldots & m_1^{d-1} \\
a_2         & m_2     & m_2^2  & \ldots & m_2^{d-1} \\
\vdots    & \vdots  & \vdots & \ddots & \vdots    \\
a_d         & m_d     & m_d^2  & \ldots & m_d^{d-1} \\
\end{vmatrix}
\end{equation}
Since the linear forms $l_i$ and $m_i$ are the equations of the line $L_{a_i}$, the polynomial $A$
belongs to the ideal of $C$. In the expansion $A=xG-zB$
the polynomial $G$ is a Vandermonde determinant
\begin{equation} \label{eqg}
G= \prod_{1 \leqslant i < j \leqslant d} (m_j-m_i) =  \prod_{1 \leqslant i < j \leqslant d} (a_i-a_j) ((a_i+a_j)z+w)
\end{equation}
We note that $G$ is nonzero: the hypothesis that the sums $a_i+a_j$ be all distinct implies
that $a_i-a_j \neq 0$ for every $i < j$. Furthermore, the zeros of $G$ in $\PP^1=\textnormal{Proj}\,(\K[z,w])$
are the points $[1:-a_i-a_j]$.

The polynomial $B$ is
\begin{equation} \label{eqh}
B= \begin{vmatrix}
a_1         & m_1     & m_1^2  & \ldots & m_1^{d-1} \\
a_2         & m_2     & m_2^2  & \ldots & m_2^{d-1} \\
\vdots    & \vdots  & \vdots & \ddots & \vdots    \\
a_d         & m_d     & m_d^2  & \ldots & m_d^{d-1} \\
\end{vmatrix}= \sum_{j=1}^d (-1)^{j-1} a_j m_1  \cdots \widehat{m_j} \cdots m_d \prod_{\begin{subarray}{c}1 \leqslant h < k \leqslant d\\ h\neq j, k \neq j \end{subarray}} (m_k-m_h)
\end{equation}

Since $m_k - m_h=  (a_h-a_k) ((a_h+a_k)z+w)$ and the initial term of $m_j$ is $y$,
the initial term of $B$ is the polynomial
\begin{equation} \label{eqlt}
\sum_{j=1}^d (-1)^{j-1} a_j y^{d-1} \prod_{\begin{subarray}{c}1 \leqslant h < k \leqslant d\\ h\neq j, k \neq j \end{subarray}}(a_h-a_k) ((a_h+a_k)z+w)=y^{d-1}P
\end{equation}
provided $P =P(z,w) \neq 0$. We will prove  not only that $P$ is not zero, but also that it has no zero in common with $G$,
that is, it does not vanish at the points $[1:-a_i-a_j]$. By symmetry, it is enough to show that
$P(1,\!-\!a_1\!-\!a_2) \neq 0$.
For this, we write $P$ as a determinant: if we set $p_i=a_i^2 z+a_i w$, then
\begin{equation} \label{eqf2}
P=
\sum_{j=1}^d (-1)^{j-1} a_j \prod_{\begin{subarray}{c}1 \leqslant h < k \leqslant d\\ h\neq j, k \neq j \end{subarray}}(p_k-p_h) =
\begin{vmatrix}
a_1     &1       & p_1        & \ldots   & p_1^{d-2} \\
a_2     &1       & p_2        & \ldots   & p_2^{d-2} \\
\vdots  &\vdots  & \vdots     & \ddots   &\vdots     \\
a_d     &1       & p_d        & \ldots   & p_d^{d-2} \\
\end{vmatrix}
\end{equation}
Now note that $p_1(1,\!-\!a_1\!-\!a_2)=p_2 (1,\!-\!a_1\!-\!a_2)=-a_1 a_2$ and
$p_i(1,\!-\!a_1\!-\!a_2)=a_i^2-(a_1+a_2)a_i$ for $j \geqslant 3$. For simplicity
we write $p_i$ in place of $p_i(1,\!-\!a_1\!-\!a_2)$.
Then
\begin{equation} \label{eqf3}
\begin{split}
P(1,\!-\!a_1\!-\!a_2)&{}=
\begin{vmatrix}
a_1     &1       & (-a_1a_2)        & \ldots   & (-a_1a_2)^{d-2} \\
a_2     &1       & (-a_1a_2)       & \ldots   & (-a_1a_2)^{d-2} \\
\vdots  &\vdots  & \vdots     & \ddots   &\vdots     \\
a_d     &1       & p_d        & \ldots   & p_d^{d-2} \\
\end{vmatrix}
=
(a_1-a_2)
\begin{vmatrix}
1       & (-a_1a_2)        & \ldots   & (-a_1a_2)^{d-2} \\
1       & p_3        & \ldots   & p_{3}^{d-2} \\
\vdots  & \vdots     & \ddots   &\vdots     \\
1       & p_d        & \ldots   & p_d^{d-2} \\
\end{vmatrix}=
\\
&{} =
(a_1-a_2) \left(\prod_{j=3}^d (p_j+a_1a_2)\right)\left(  \prod_{3 \leqslant h < k \leqslant d} (p_k-p_h) \right)= \\
&{}=
(a_1-a_2) \left(\prod_{j=3}^d (a_j-a_1)(a_j-a_2) \right)\left(  \prod_{3 \leqslant h < k \leqslant d} (a_k\!-\!a_h)(a_h+a_k\!-\!a_1\!-\!a_2) \right).
\end{split}
\end{equation}
This shows $P(1,\!-\!a_1\!-\!a_2) \neq 0$ because of the assumption that sums $a_i+a_j$ be all distinct.

To finish, we let $F=zP$. Then $\text{in}_{\omega}(A)=xG-y^{d-1}F$, and $F$ and $G$ are homogeneous forms in $\K[z,w]$, $\deg G= \binom{d}{2}$, $\deg F= \binom{d-1}{2}+1$,
and $F$ and $G$ have no common zero in $\PP^1$. This concludes the proof.

We would like to make some remarks on the polynomial $P$. It is divisible by $z^{d-2}$. Indeed,
$P$ is a form in $z$ and $w$ of degree $ \binom{d-1}{2}$
with coefficients in $\K[a_1, \ldots,a_d]$. It is divisible by the Vandermonde determinant $V(a_1,\ldots,a_d)$
because it is antisymmetric in the $a_i$'s. Furthermore, the coefficient of $z^m w^n$ in $P$ is an antisymmetric polynomial
of degree $2m+n+1$ in the $a_i$'s: in order for it to be nonzero, it is necessary that
$2 m + n +1\geqslant \deg V= \binom{d}{2}$. Since $ m+n = \binom{d-1}{2}$, we deduce
$ m \geqslant d-2$. This means that $P$ is divisible by $z^{d-2}$, and the coefficient of $z^{d-2}w^{\binom{d-2}{2}}$
is  $ V(a_1, \ldots, a_d)$ times a constant $-c_d$ that depends on $d$ but not on the $a_i$'s.  One can compute
$c_2=1$ and
 $   c_d = {\ds \sum_{k=1}^{d-2}(-1)^{k+1}} \binom{d-1-k}{k}c_{d-k}$ for $d \geqslant 3$.
The first few values are $c_3=1$, $c_4=2$, $c_5=5$, $c_6=14$, $c_7=42$, and in general
$c_d$ is the $(d-2)^{th}$ Catalan number:
$
c_d = C(d-2)= \frac{1}{d-1}\binom{2d-4}{d-2}$.
\end{proof}

\br
We can give a geometric interpretation of the condition that the sums
$a_i+a_j$ be all distinct. In the family constructed in the proof of the theorem,
the union of the two lines $L_{a_i}$ and $L_{a_j}$ specializes to the planar double line
$x=y^2=0$ plus the embedded point  $x=y= (a_i+a_j)z+w=0$.
Thus the condition means that these embedded points are all distinct.
\er

\bex
If the condition that the $a_i+a_j$ be all distinct is not satisfied, we expect the limit to acquire embedded
points. The reason is that in this case the proof of Theorem \ref{main} shows
that the polynomials $F$ and $G$ have a common zero, and so $xG-y^{d-1}F$ is no longer irreducible.
For a specific example, we take $d=4$ and $a_1=0$, $a_2=1$, $a_3=2$ and $a_4=3$ (in characteristic $\neq 2,3$),
so that
$a_1+a_4=a_2+a_3=3$.  In this case,
\[
\begin{split}
A = xG-zB ={}&  x \big(12(z+w)(2z+w)(3z+w)^2(4z+w)(5z+w)\big)- \\ &z \big(12 y(3z+w)(2 y^{2} z^{2}-30 y z^{3}+148 z^{4}-15 y z^{2} w+195 z^{3} w-y z w^{2}+85 z^{2} w^{2}+15 z w^{3}+w^{4})\big)
\end{split}
\]
and
\[
\begin{split}
\text{in}_{(4,2,1,1)} (A)&{} = xG-y^3F = x \big(12(z+w)(2z+w)(3z+w)^2(4z+w)(5z+w)\big)- y^3 \big(24 z^3(3z+w)\big) ={}\\
&{} = 12(3z+w)\big(xG_1 - y^3 F_1\big).
\end{split}
\]
The initial ideal of $I_C$, according to \textit{Macaulay2}, is
\[
\text{in}_{(4,2,1,1)}(I_C) = \big\langle x^{2},6 x y z^{2}+2 x y z w,y^{4},x y^{2} w,6 x
      y^{2} z,6 x y z w^{2}+2 x y w^{3},xF_1 - y^3 G_1,6 x y^{3}\big\rangle
\]
whose saturation is
\[
\big\langle x^{2}, x y(3z+w),x y^{2},y^{4},xG_1 - y^3 F_1\big\rangle.
\]
Therefore the limit $C_0$ consists of the unique $4$ structure supported on $L$
contained in the surface $xG_1-y^3F_1$, which by Proposition \ref{extr_mls} is an extremal curve
of genus $-2$, plus an embedded point (of equation $3z+w=0$ on $L$).
\eex

\bc
Let $C$ be an effective divisor of bidegree $(d,0)$ on a smooth quadric surface. Then every curve in the biliaison class
of $C$ is in the connected component of the extremal curves in its Hilbert scheme. In particular,
a curve $D$ on a smooth quadric surface is in the connected component of the extremal curves in its Hilbert scheme.
\ec

\begin{proof}
The case $0 \leqslant d \leqslant 2$ has been proven
by Hartshorne \cite{hconn}. To prove the statement for $d\geqslant 3$, 
we may assume that the quadric surface has equation as in  Theorem \ref{main}; the claim then
follows from \ref{main} and \cite[Theorem 2.3]{footnote}
because $C$ is a minimal curve. Finally, every curve $D$ on a smooth quadric surface $Q$ is obtained by an ascending biliaison
on $Q$ from a curve of bidegree $(d,0)$ (or $(0,d)$, which reduces to the previous case by a change of variables).
\end{proof}

\end{document}